\documentclass[12pt]{amsart}
\usepackage[latin9]{inputenc}
\usepackage[a4paper]{geometry}
\usepackage{textcomp}
\usepackage{enumitem}
\usepackage{amsbsy}
\usepackage{amstext}
\usepackage{amsthm}
\usepackage{amssymb}

\makeatletter

\newcommand{\noun}[1]{\textsc{#1}}
\DeclareTextSymbolDefault{\textquotedbl}{T1}

\numberwithin{equation}{section}
\numberwithin{figure}{section}
\theoremstyle{remark}
\newtheorem*{acknowledgement*}{\protect\acknowledgementname}
\theoremstyle{definition}
\newtheorem*{defn*}{\protect\definitionname}
\theoremstyle{plain}
\newtheorem*{thm*}{\protect\theoremname}

\@ifundefined{date}{}{\date{}}
\usepackage[nobysame,abbrev,bibtex-style]{amsrefs}
\usepackage{undertilde}
\usepackage{amsfonts}
\usepackage{bbm}
\setlist[enumerate]{leftmargin=*,widest=0}
\setitemize[0]{leftmargin=10pt,itemindent=0pt}

\theoremstyle{definition}
\newtheorem*{thma1*}{Theorem A1}\newtheorem*{thma2*}{Theorem A2}\newtheorem*{thmb*}{Theorem B}\newtheorem*{thmc*}{Theorem C}\def\be{\begin{equation}}
\def\ee{\end{equation}}
\def\be*{\begin{equation*}}
\def\ee*{\end{equation*}}

\def\bbr{\mathbb{R}}
\def\bbh{\mathbb{H}}
\def\la{\lambda}

\numberwithin{equation}{section}

\def\PGL{\text{\rm PGL}}
\def\SL{\text{\rm SL}}
\def\GL{\text{\rm GL}}

\def\Ga{\Gamma}
\def\bbr{\mathbb{F}}

\def\bbz{\mathbb{Z}}

\def\la{\lambda}

  \def\ee{\end{equation}}

\newcommand{\one}{\mathbbm{1}}

\newcommand{\Mod}[1]{\,\left(\textup{mod}\;#1\right)}

\theoremstyle{definition} 
\def\bbq{\mathbb{Q}}

\def\bbr{\mathbb{R}}
\def\bba{\mathbb{A}}

\def\bbh{\mathbb{H}}
\def\be{\begin{equation}}
\def\ee{\end{equation}}

\def\vare{\varepsilon}

\theoremstyle{remark}  

\makeatother

\usepackage{listings}
\providecommand{\acknowledgementname}{Acknowledgement}
\providecommand{\definitionname}{Definition}
\providecommand{\theoremname}{Theorem}

\begin{document}
\title[From Ramanujan Graphs to Ramanujan Complexes]{From Ramanujan Graphs\\
to Ramanujan Complexes}
\author{Alexander Lubotzky and Ori Parzanchevski}
\begin{abstract}
Ramanujan graphs are graphs whose spectrum is bounded optimally. Such
graphs have found numerous applications in combinatorics and computer
science. In recent years, a high dimensional theory has emerged. In
this paper these developments are surveyed. After explaining their
connection to the Ramanujan conjecture we will present some old and
new results with an emphasis on random walks on these discrete objects
and on the Euclidean spheres. The latter lead to \textquotedbl golden
gates\textquotedbl{} which are of importance in quantum computation.
\end{abstract}

\maketitle
\setcounter{section}{-1}

\section{Introduction}

Let $X$ be a finite connected $k$-regular graph and $A$ its adjacency
matrix. The graph $X$ is called Ramanujan graph if every eigenvalue
$\la$ of $A$ satisfies either $|\la|=k$ or $|\la|\le2\sqrt{k-1}$.
This term was coined in \cite{LPS88}.

While Ramanujan had an interest in combinatorics (the partition function
etc.), it does not seem as though he has had a special interest in
graph theory. So why are these graphs named after him? This will be
explained in \S1. The explanation will suggest what should be the
definition of Ramanujan graph, for a general graph, not necessarily
$k$-regular. Moreover, it will suggest the definition for directed
graphs (digraphs) and even high dimensional simplicial complexes (the
so called Ramanujan complexes), as will be explained in \S 2 and
in \S3.

Ramanujan graphs have found plenty of applications in computer science
and pure mathematics. Most of them have to do with the fact that they
provide optimal expanders (see \cite{HLW06,Lub10,Lub12} and the references
therein). Lately, Ramanujan complexes and high dimensional expanders
have also started to be a popular subject of research (cf.\ \cite{Lubotzky2013,lubotzky2017high}
and the references therein).

Here we concentrate on describing their aspects which truly use the
full power of being Ramanujan, and not merely expansion: In §4, we
will describe random walks on Ramanujan graphs and complexes and in
§5, ``golden gates'', which is a new fascinating application of
them to quantum computation.
\begin{acknowledgement*}
The authors acknowledge with gratitude a support by the ERC and the
NSF (A.L.), and the ISF (O.P.).
\end{acknowledgement*}

\section{Why Ramanujan?}

Let $X$ be a finite connected $k$-regular graph, $k\ge3$, with
$n$ vertices, and $A$ its adjacency $n\times n$ matrix. Being symmetric,
all its eigenvalues $\lambda$ are real and it is easy to see that
$|\lambda|\le k$, $k$ is always an eigenvalue, and $-k$ is an eigenvalue
if and only if $X$ is bi-partite. The graph $X$ is called \emph{Ramanujan
graph} if every eigenvalue $\lambda$ satisfies either $|\lambda|=k$
or $|\lambda|\le2\sqrt{k-1}$. The bound $2\sqrt{k-1}$ is significant:
by Alon-Boppana Theorem, (cf.\ \cite[Prop.\ 4.2]{LPS88}) this is
the best possible bound one can hope for, for an infinite family of
$k$-regular graphs. The real reason behind it is as follows: The
universal cover of $X$ (in the sense of algebraic topology) is $\tilde{X}=T_{k}$
- the infinite $k$-regular tree. An old result of Kesten asserts
that the spectrum of the adjacency operator acting on $L^{2}(T_{k})$
is the interval $[-2\sqrt{k-1},2\sqrt{k-1}]$. So, being Ramanujan
means for $X$, that all its non-trivial eigenvalues are in the spectrum
of its universal cover $\tilde{X}$.

Ramanujan graphs are optimal expanders from spectral point of view.
Recall that a finite $k$-regular graph $X$ is called $\vare$-expander
if $h(X)\ge\vare$ when $h(X)$ is the Cheeger constant of $X$, namely
\[
h(X)=\min\left\{ \frac{|E(Y,\bar{Y})|}{|Y|}\,\middle|\,{Y\subseteq X,\atop 0<|Y|\le\frac{|X|}{2}}\right\} 
\]
when $E(Y,\bar{Y})$ is the set of edges between $Y$ and its complement.

Now if we denote $\lambda_{1}(X)=\max\{\lambda\,|\,\lambda\neq k,\lambda\text{ e.v.\  of }A\}$,
then
\[
\frac{h(X)^{2}}{2k}\le k-\lambda_{1}(X)\le2h(X)\hbox{\ (cf.\ \cite[\S4.2]{Lub10})}.
\]
So, Ramanujan graphs are expanders. Expander graphs are of great importance
in combinatorics and computer science (cf.\ \cite{HLW06} and the
references therein) and also in pure mathematics (cf.\ \cite{Lub12}).
Expander graphs serve as basic building blocks in various network
constructions, in many algorithms and so on. The bound on their eigenvalues
ensures that the random walk on such graphs converges quickly to the
uniform distribution and on Ramanujan graphs this happens in the fastest
possible way (see §4). This is one more reason that makes them so
useful.

The existence of Ramanujan graphs is by no means a trivial issue:
While it is known that random $k$-regular graphs are expanders, it
is not known if they are Ramanujan. First examples of infinite families
of such graphs were given by explicit construction in \cite{LPS88}
and \cite{margulis1988explicit} for $k=q+1$, $q$ prime. In \cite{marcus2013interlacing},
it is shown, by a non constructive method, that for every $k\ge3$
there exist infinitely many $k$-regular bi-partite Ramanujan graphs.

Why are Ramanujan graphs named after Ramanujan? As far as we know
Ramanujan had no special interest in graph theory. Let us explain
the reason for this name which was coined in \cite{LPS88}.

Observe the following power series 
\[
\Delta(q)=q\prod\limits _{n\ge1}(1-q^{n})^{24}=\Sigma\tau(n)q^{n}=q-24q+252q^{3}+\dots
\]

The coefficients $\tau(n)$ define the so called Ramanujan tau function.
Ramanujan conjectured that $|\tau(p)|\le2p^{\frac{11}{2}}$ for every
prime $p$. The importance of $\Delta$ comes from the fact that if
we write $q=e^{2\pi iz}$ then $\Delta(z)$ is a cusp form of weight
12 on the upper half plane ${\bbh}=\{z=x+iy\,|\,x,y\in\bbr,\;\;y>0\}$
with respect to the modular group $\Gamma=\SL_{2}(\bbz)$ acting on
$\bbh$ by Möbius transformation ${a\:b \choose c\:d}(z)=\frac{az+b}{cz+d}$.
Now if $\Gamma_{0}(N)=\{{a\:b \choose c\:d}\in\Gamma\,\Big|\,c\equiv0\!\mod N\}$
we denote $S_{k}(N)$ (or more generally $S_{k}(N,w)$ for a Dirichlet
character $w$ of $\bbz/N\bbz$) the space of cusp forms on $\bbh$
w.r.t. $\Ga_{0}(N)$ (and $w$). The Hecke operators $T_{p}$ ($p$
prime, $(p,N)=1$), act, and commute, on each $S_{k}(N,w)$, and their
common eigenfunctions are the Hecke eigenforms. Now, $S_{12}(\Ga=\Ga_{0}(1))$
is one dimensional and so $\Delta(z)$ above is such a Hecke eigenform.
Moreover, $\tau(p)$ above is equal to the eigenvalue of $T_{p}$
acting on $S_{12}(\Gamma)$. A natural and far reaching generalization
of the Ramanujan conjecture mentioned above on the size of $\tau(p)$
is the so called Ramanujan-Peterson (RP) conjecture: for every Hecke
eigenform $f$ in $S_{k}(N,w)$, the eigenvalues $\lambda_{p}$ of
$T_{p}$, $(p,N)=1$, satisfy $|\lambda_{p}|\le2p^{\frac{k-1}{2}}$.
The reader is referred to \cite{RogawskiModularformsRamanujan} for
a concise and clear explanation of all these notions.

The modern approach to automorphic functions via representation theory
brought in another point of view on the Ramanujan-Peterson Conjecture.
Satake \cite{satake1966spherical} showed that the RP conjecture is
equivalent to the assertion: Let $\bba=\bba_{\bbq}$ be the ring of
adeles of $\bbq$, and $\pi$ an irreducible cuspidal $\GL_{2}$-representation
in $L^{2}(\GL_{2}(\bba)/\GL_{2}(\bbq))$, such that its component
at infinity $\pi_{\infty}$ is square integrable, then for every prime
$p$ the local factor at the $p$-component $\pi_{p}$ is a tempered
representation. See \cite{RogawskiModularformsRamanujan} for exact
definitions. Here we only mention that a representation of a (simple)
$p$-adic or real Lie group $G$ is tempered if it is weakly contained
in $L^{2}(G)$. The RP conjecture was proved by Deligne (for the special
representations that are relevant to the Ramanujan graphs, the RP
conjecture was actually proven earlier by Eichler). The representation
theoretic formulation suggests vast generalizations to other algebraic
groups.

Let us look at the $p$-adic group $G=\PGL_{2}(\bbq_{p})$. The Bruhat-Tits
building associated with $G$ is, in this special case, the $(p+1)$-regular
tree $T=T_{p+1}$ which can be identified as $T=G/K$ when $K$ is
a maximal compact subgroup of $G$. If $\Gamma$ is a discrete cocompact
subgroup of $G$, then $X=\Gamma\setminus T=\Gamma\setminus G/K$
is a finite $(p+1)$-regular graph. One can show (see \cite{Lub10})
that $X$ is a Ramanujan graph if and only if every infinite dimensional
$K$-spherical $G$-sub-representation of $L^{2}(\Gamma\setminus G)$
is tempered. Deligne theorem, combined with the so called Jacquet-Langlands
correspondence, enables the construction of such arithmetic subgroups
$\Ga$ for which the temperedness condition is satisfied and hence
Ramanujan graphs are obtained. This was the method of \cite{LPS88}
and \cite{margulis1988explicit}. Let us mention that for every $k$,
if $G$ is the full automorphism group of $T_{k}$ and $\Gamma$ a
discrete cocompact subgroup of $G$, then $X=\Ga\setminus T_{k}$
is $k$-regular Ramanujan graph if and only if the same temperedness
condition is satisfied: in other words every non-trivial eigenvalue
of $X=\Ga\setminus T_{k}$ is coming from the spectrum of $T_{k}$
if and only if every non-trivial spherical subrepresentation of $L^{2}(\Ga\setminus G)$
is coming from $L^{2}(G)$. This illustrates the connection between
the notion of Ramanujan graph and the Ramanujan conjecture.

As mentioned above, the Ramanujan-Peterson conjecture was generalized
to other groups, and some of its generalizations to $\GL_{d}$ (instead
of only $\GL_{2}$) led to higher dimensional versions of Ramanujan
graphs, the so called Ramanujan complexes. We will see more on it
in §3.

Finally, another interesting hint to a connection with number theory:
Ihara defined the notion of Zeta function of a $k$-regular graph
$X$, and Sunada observed that $X$ is Ramanujan if and only if this
Zeta function satisfies ``the Riemann hypothesis''. We refer the
reader to the survey \cite{LiRamanujanconjectureits} for more details.

\section{General graphs and digraphs}

The first paragraph of §1 suggests what should be the general definition
of Ramanujan graphs. This was carried out for the first time in Greenberg
Thesis (\cite{greenberg1995spectrum}, which is unfortunately not
published and available only in Hebrew), and was vastly generalized
in \cite{grigorchuk1999asymptotic}.

Here is the main point. Let $X$ be any finite connected graph and
$\widetilde{X}$ its universal cover. Let $A_{\widetilde{X}}$ be
the adjacency operator acting on $L^{2}(\widetilde{X})$ by $A_{\widetilde{X}}\left(f\right)\left(x\right)=\sum_{x'\sim x}f\left(x'\right)$
where $x'$ runs over the neighbors of $x$ in $\widetilde{X}$. Now,
it is shown in \cite{greenberg1995spectrum} that there exists a positive
real number $\kappa$ depending only on $\widetilde{X}$, such that
if $Y$ is a finite graph covered by $\widetilde{X}$, then $\kappa$
is the largest (Perron-Frobenius) eigenvalue of the adjacency matrix
$A_{Y}$ of $Y$. When $X$ is $k$-regular $\kappa=k$, and when
$X$ is a bipartite $\left(k_{1},k_{2}\right)$-biregular, $\kappa=\sqrt{k_{1}k_{2}}$.
\begin{defn*}[\cite{greenberg1995spectrum}]
The graph $X$ is called Ramanujan if every eigenvalue $\lambda$
of $A_{X}$ satisfies either $|\lambda|=\kappa$ or $\lambda\in\mathrm{Spec}(A{}_{L^{2}(\widetilde{X})})$.
\end{defn*}
This recovers the classical definition of Ramanujan graphs for $k$-regular
graphs since $\mathrm{Spec}(A|_{L^{2}(T_{k})})=\left[-2\sqrt{k-1},2\sqrt{k-1}\right]$.
For bipartite $(k_{1},k_{2})$-biregular graphs $X$ with $k_{1}\leq k_{2}$,
the universal cover is the $(k_{1},k_{2})$-biregular tree $T_{k_{1},k_{2}}$
and 
\begin{multline*}
\mathrm{Spec}(A\big|_{L^{2}(T_{k_{1},k_{2}})})={\textstyle \left[-\sqrt{k_{2}-1}-\sqrt{k_{1}-1},-\sqrt{k_{2}-1}+\sqrt{k_{1}-1}\right]}\\
{\textstyle \cup\left\{ 0\right\} \cup\left[\sqrt{k_{2}-1}-\sqrt{k_{1}-1},\sqrt{k_{2}-1}+\sqrt{k_{1}-1}\right].}
\end{multline*}
It is known that for every $3\leq k\in\mathbb{N}$, there exist infinitely
many $k$-regular Ramanujan graphs (explicit constructions for every
$k=p^{e}+1$, $p$ prime \cite{morgenstern1994existence}, and non
explicit for every $k$ \cite{marcus2013interlacing}). But for $(k_{1},k_{2})$-biregular,
it is known only for special values:
\begin{thm*}[\cite{ballantine2015explicit,Evra2018RamanujancomplexesGolden}]
Let $p$ be a prime, $k_{1}=p+1$ and $k_{2}=p^{3}+1$, then there
exist infinitely many bipartite $(k_{1},k_{2})$-biregular Ramanujan
graphs.
\end{thm*}
In \cite{ballantine2015explicit} existence was shown as the quotients
of the bi-regular tree associated with a rank one simple $p$-adic
Lie group. Explicit constructions (in the sense of computer science)
are given for $p\equiv3\negthickspace\mod(4)$ in \cite{Evra2018RamanujancomplexesGolden}.

Let us mention that \cite{marcus2013interlacing} gives existence
of ``weak-Ramanujan'' $(k_{1},k_{2})$-biregular graphs in the following
sense: every eigenvalue $\lambda$ is either $|\lambda|=\kappa=\sqrt{k_{1}k_{2}}$
or $|\lambda|\leq\sqrt{k_{2}-1}+\sqrt{k_{1}-1}$.

In \cite{Lubotzky1998Noteveryuniform} it was shown that there exist
finite graphs $X$ for which $\widetilde{X}$ does not cover any Ramanujan
graph. This was put in a more general framework in \cite{friedman2003relative}.

Turning to digraphs (directed graph), denote by $A=A_{\mathcal{D}}$
the adjacency matrix of the digraph $\mathcal{D}$, namely $A_{v,w}=1$
if $v\rightarrow w$ in $\mathcal{D}$ and $A_{v,w}=0$ otherwise.
We say that $\mathcal{D}$ is $k$-regular if every vertex has $k$
incoming edges, and $k$ outgoing ones. The notion of Ramanujan digraphs
(directed graphs) was considered only quite recently \cite{Parzanchevski2018SuperGoldenGates,Lubetzky2017RandomWalks,Parzanchevski2018RamanujanGraphsDigraphs}.
A main reason for this is that the adjacency matrix of a digraph can
be non-normal, in which case its spectrum reveals much less information
on the graph.
\begin{defn*}
A $k$-regular digraph is a Ramanujan digraph if every eigenvalue
of $A_{X}$ satisfies either $\left|\lambda\right|=k$ or $\left|\lambda\right|\leq\sqrt{k}$.
\end{defn*}
Here the trivial eigenvalues can be $e^{2\pi i/m}k$ for any $m\in\mathbb{N}$,
indicating that the digraph is \emph{$m$-periodic}: its vertices
can be partitioned into $m$ sets $V_{0},\ldots,V_{m-1}$, with every
edge starting in $V_{j}$ terminating in $V_{\left(j+1\!\mod m\right)}$.
Once again, the non-trivial spectrum agrees with the ``directed universal
cover'' $T_{k}^{\rightleftharpoons}$, which is the $2k$-regular
tree, directed to have constant in-degree and out-degree $k$. Indeed,
$\mathrm{Spec}(A\big|_{L^{2}(T_{k}^{\rightleftharpoons})})=\{z\in\mathbb{C}|\left|z\right|\leq\sqrt{k}\}$
by \cite{Harpe1993spectrumsumgenerators}.

A general example of a Ramanujan digraph arises from Hashimoto's approach
to Ihara's zeta function \cite{hashimoto1989zeta}. Given a $\left(k+1\right)$-regular
(undirected) graph $X$, define the $k$-regular digraph $D_{X}$,
whose vertices correspond to directed edges in $X$, and whose edges
correspond to non-backtracking steps in $X$. Namely, $e\rightarrow e'$
in $D_{X}$ iff $e,e'$ form a non-backtracking path in $X$. Hashimoto's
work shows that $D_{X}$ is a Ramanujan digraph if and only if $X$
is a Ramanujan graph.

It is interesting to note that the Alon-Boppana theorem fails for
digraphs: the \emph{De-Bruijn }digraphs (cf.\ \cite[§3.4]{Parzanchevski2018RamanujanGraphsDigraphs})
are $k$-regular digraphs, of arbitrarily large size, whose non-trivial
spectrum consists entirely of zeros! However, these graphs have non-normal
adjacency matrices. It turns out that normality, and even ``almost-normality''
recovers an Alon-Boppana bound, for which Ramanujan digraphs are again
optimal. We say that a family of digraphs is \emph{almost-normal }if
the adjacency matrices of its members are unitarily equivalent to
block-diagonal matrices with blocks of globally bounded size.
\begin{thm*}[\cite{Parzanchevski2018RamanujanGraphsDigraphs}]
The smallest upper bound for the non-trivial spectrum of an infinite
almost-normal family of $k$-regular, $m$-periodic digraphs, is $\sqrt{k}$.
\end{thm*}
It turns out that almost-normality appears naturally in the context
of digraphs which arise from Ramanujan graphs and complexes (see §\ref{sec:Ramanujan-complexes}),
and that it serves as a substitute for normality in the spectral analysis
of these digraphs.

\section{\label{sec:Ramanujan-complexes}Ramanujan complexes}

Combinatorial graphs are one-dimensional simplicial complexes, and
it is natural to ask for analogues of expanders and Ramanujan graphs
in higher dimension. Here even the definition is not straightforward,
as there is no clear counterpart to the $k$-regular tree $T_{k}$
in general dimension. The explicit construction of Ramanujan graphs
suggests one answer: since for $k=p+1$ the tree $T_{k}$ arose as
the Bruhat-Tits building of $G=PGL_{2}\left(\mathbb{Q}_{p}\right)$,
one can replace it with the Bruhat-Tits building $\mathcal{B}=\mathcal{B}\left(G\right)$
of $G=PGL_{d+1}\left(\mathbb{Q}_{p}\right)$, which is an infinite,
contractible, $d$-dimensional complex. This is indeed the approach
taken in \cite{li2004ramanujan,Lubotzky2005a}, except for the replacement
of $\mathbb{Q}_{p}$ by $\mathbb{F}_{p}\left(\left(t\right)\right)$
- the reason being that the Ramanujan conjecture for $PGL_{d}$ over
$\mathbb{Q}$ is still open for $d\geq3$, whereas for $PGL_{d}$
over $\mathbb{F}_{p}\left(t\right)$ it was proved by Lafforgue in
\cite{lafforgue2002chtoucas}. A more general approach is to look
at any non-archimedean local field $F$, and $G=\utilde{G}\left(F\right)$,
where $\utilde{G}$ is a simple $F$-algebraic group. Bruhat-Tits
theory associates with $G$ a building $\mathcal{B}$ (the so called
Bruhat-Tits building) which is a contractible simplicial complex of
dimension $d$ equal to the $F$-rank of $\utilde{G}$. The group
$G$ acts on $\mathcal{B}$, transitively on the $d$-cells. Every
torsion-free discrete cocompact subgroup $\Gamma$ of $G$ gives rise
to a finite complex $X=\Gamma\backslash\mathcal{B}$, which can then
be compared to its universal cover $\widetilde{X}=\mathcal{B}$.

For this comparison, one should decide which adjacency operator should
one look at, as the standard adjacency relation between vertices depends
only on the 1-skeleton of the complex, and does not capture the high-dimensional
structure. One can ask, for example, about operators such as the discrete
$j$-dimensional Laplacian, which acts on cells in dimension $j$
and detects the presence of real $j$-th cohomology. We take an inclusive
approach: we call an operator $T$ on (a subset of) the cells of the
building $\mathcal{B}\left(G\right)$ \emph{geometric} if it commutes
with the action of $G$. If $X$ is a finite quotient of $\mathcal{B}$,
this implies that $T$ descends to a well-defined operator $T|_{X}$
on $X$, and we define:
\begin{defn*}
Let $F$ be a nonarchimedean local field, $\mathcal{B}$ the Bruhat-Tits
building associated with $PGL_{d+1}\left(F\right)$, and $X$ a quotient
of $\mathcal{B}$.
\begin{enumerate}
\item For a geometric operator $T$, an eigenvalue of $T|_{X}$ is \emph{trivial
}if the associated eigenfunction on $X$ lifts to a $PSL_{d+1}\left(F\right)$-invariant
function on $\mathcal{B}$.
\item The complex $X$ is a \emph{Ramanujan complex }if for every geometric
operator $T$ on $\mathcal{B}$, the nontrivial spectrum of $T|_{X}$
is contained in the $L^{2}$-spectrum of $T$ on $\mathcal{B}$.
\end{enumerate}
\end{defn*}
The definition generalizes to other groups than $PGL_{d}$, once we
understand which are the trivial eigenfunctions - see \cite{Lubetzky2017RandomWalks}
for the case of simple algebraic groups, and \cite{first2016ramanujan}
for a more general one.

We remark that the original definition of Ramanujan complexes in \cite{li2004ramanujan,Lubotzky2005a}
only requires (2) for geometric operators on the vertices of $\mathcal{B}$.
However, all the known constructions of Ramanujan complexes \cite{li2004ramanujan,Lubotzky2005b,sarveniazi2007explicit,first2016ramanujan,Evra2018RamanujancomplexesGolden}
satisfy the stronger definition!

As in the case of graphs, the Ramanujan property can be related to
representation theory: The Iwahori group of $G$ is the pointwise
stabilizer of a cell of maximal dimension in $\mathcal{B}$, and the
complex $X=\Gamma\backslash\mathcal{B}$ is Ramanujan if and only
if every infinite dimensional, Iwahori-spherical, irreducible $G$-sub-representation
of $L^{2}(\Gamma\setminus G)$ is tempered \cite{kamber2016lp,first2016ramanujan,Lubetzky2017RandomWalks}.

\section{Random walks}

A highly useful property of expanders is that random walks on them
converge rapidly to the stationary distribution: Let $X$ be a non-bipartite
$k$-regular graph, $\left\{ v_{t}\right\} _{t=0}^{\infty}$ a simple
random walk (SRW) process on $X$, and $P_{X}^{t}:v\mapsto Prob\left[v_{t}=v\right]$
the distribution of the walk at time~$t$. It is a standard exercise
to show that $\left\Vert P_{X}^{t}-\mathbf{u}\right\Vert _{2}$, the
$L^{2}$-distance of $P_{X}^{t}$ from the uniform distribution, is
bounded by $\left(\lambda/k\right)^{t}$, where $\lambda$ is the
largest non-trivial eigenvalue of $A_{X}$ (in absolute value).

It turns out, however, that Ramanujan graphs are optimally mixing
not only in $L^{2}$-norm, but also in $L^{p}$ for all $1\leq p\leq\infty$.
Furthermore, they manifest a \emph{cutoff} phenomena: the $L^{p}$-distance
$\left\Vert P_{X}^{t}-\mathbf{u}\right\Vert _{p}$ drops abruptly
from being near maximal to being near zero, over a short interval
of time called the \emph{cutoff window}. We focus on $L^{1}$, the
\emph{total-variation} norm, which is hardest to bound, and the most
useful for many purposes (see \cite{Lovasz1996Randomwalksgraphs}).
\begin{thm*}[\cite{lubetzky2016cutoff}]
Let $X$ be a $k$-regular Ramanujan graph on $n$ vertices.
\begin{enumerate}
\item The SRW on $X$ has $L^{1}$-cutoff at time $\frac{k}{k-2}\log_{k-1}n$.
\item The non-backtracking random walk (NBRW) on $X$ has $L^{1}$-cutoff
at time $\log_{k-1}n$.
\end{enumerate}
\end{thm*}
Notice that the location of cutoff for NBRW is optimal: a non-backtracking
walker on a $k$-regular graph sees at most $k-1$ new vertices at
every step, with the exception of the first one. Thus, a walk of length
$\left(1-\delta\right)\log_{k-1}n$ can reach only a small fraction
of the graph for $\delta>0$ (and even for $\delta=1/\log\log n$),
resulting in $L^{1}$-distance $1-o\left(1\right)$ from equilibrium.
In a similar manner one can show that the first bound is optimal when
taking into account the hindrance caused by backtracking.

In \cite{lubetzky2016cutoff}, the authors first prove the bound for
NBRW on $X$, and then show that it implies the bound for SRW. Let
us give a glimpse of how the bound for NBRW is proved. Recall the
digraph $D_{X}$ from §2: this is a \textbf{$\boldsymbol{\left(k\!-\!1\right)}$}-regular
Ramanujan digraph, and by its construction SRW on $D_{X}$ is equivalent
to NBRW on $X$. If the adjacency matrix $A_{D_{X}}$ was symmetric,
or even normal, then we would have $\left\Vert P_{D_{X}}^{t}-\mathbf{u}\right\Vert _{2}\leq\left(k-1\right)^{-t/2}$
as for undirected expanders, and a standard $L^{2}$ to $L^{1}$ bound
would then give the desired result. However, $A_{D_{X}}$ is not normal
when $k\geq3$. The main step in \cite{lubetzky2016cutoff} is to
show that $D_{X}$ is $2$-normal, namely, $A_{D_{X}}$ is unitarily
equivalent to a block-diagonal matrix with blocks of size $2\times2$.
This is then shown to imply the bound $\left\Vert P_{D_{X}}^{t}-\mathbf{u}\right\Vert _{2}\leq\left(t+1\right)\left(k-1\right)^{-(t+1)/2}$,
which only differs by a logarithmic factor, and suffices to prove
cutoff.

Let us stress that the work of Lubetzky and Peres \cite{lubetzky2016cutoff}
uses the full strength of the Ramanujan property to deduce the cutoff
phenomenon. It is still a widely open conjecture of Peres that such
phenomena happens in all transitive expander graphs. It is known that
it is not always the case for general expanders \cite{lubetzky2011explicit}.

In Ramanujan complexes of higher dimension, it turns out that the
digraph $D_{X}$ induced by NBRW is not a Ramanujan digraph anymore.
However, it is shown in \cite{Lubetzky2017RandomWalks} that other
operators on the cells of these complexes do induce Ramanujan digraphs.
The crucial property is that these operators should describe \emph{collision-free
}walks on the building: this means that all the paths which descend
from a fixed starting cell never meet one another (for example, non-backtracking
walk on a tree have this property). It is shown in \cite{Lubetzky2017RandomWalks}
that if a geometric operator induces a collision-free walk on $\mathcal{B}$,
and $X$ is a Ramanujan quotient of $\mathcal{B}$, then the digraph
which represents the walk by $T$ on $X$ is a Ramanujan digraph.
Furthermore, it is shown that the digraphs which arise from quotients
of a fixed building are almost-normal, which leading again to cutoff
at the optimal time:
\begin{thm*}[\cite{Lubetzky2017RandomWalks}]
Let $T$ be a geometric, $k$-regular, collision-free operator on
$\mathcal{B}$, the Bruhat-Tits building of a simple $p$-adic group
$G$. Then the walk induced by $T$ on a Ramanujan complex $X=\Gamma\backslash\mathcal{B}$
has $L^{1}$-cutoff at time $\log_{k}|X|$.
\end{thm*}
In addition, it is shown in \cite{Lubetzky2017RandomWalks} that such
walks do exist: for $G=PGL_{d+1}\left(F\right)$, a collision-free
walk on $j$-cells is exhibited for each $1\leq j\leq d$, the so
called \emph{geodesic $j$-flow}. For example, geodesic $1$-flow
goes from a (colored) edge $v\rightarrow w$ to $w\rightarrow u$
if the cell $\left\{ v,w,u\right\} $ does \emph{not }belong to the
complex. The situation when $j=0$ is different: due to commutativity
of the Hecke algebra, no geometric operator on vertices induces a
Ramanujan digraph (see \cite[Rem.\ 3.5(b)]{Parzanchevski2018RamanujanGraphsDigraphs}).
However, it is shown in \cite{Chapman2019CutoffRamanujancomplexes}
that by combining the optimal cutoff result for the $j$-flow operators
in all dimensions, it is possible to recover cutoff for SRW on vertices:
\begin{thm*}[\cite{Chapman2019CutoffRamanujancomplexes}]
SRW on the vertices of Ramanujan complexes associated with $PGL_{d}\left(F\right)$
exhibit $L^{1}$-cutoff.
\end{thm*}
Once again, the proof requires the strength of the Ramanujan property,
and not merely expansion. Moreover, it needs the full high dimensional
structure of X, even when we study the SRW only on the vertices.

Finally, we mention that in \cite{golubev2017cutoff} a different
direction is taken: replacing $PGL_{2}\left(\mathbb{Q}_{p}\right)$
with $PGL_{2}\left(\mathbb{R}\right)$, the authors suggest the notion
of \emph{Ramanujan surface}s, which are hyperbolic Riemann surfaces
which spectrally behave like their universal cover, the hyperbolic
plane. It is then shown that a discrete random walk with constant-length
steps on these surfaces exhibits $L^{1}$-cutoff.

\section{Golden Gates}

Recently, Ramanujan graphs and complexes have found a surprising application
to the theory of quantum computation. In classical computation, one
decomposes any function into basic logical gates such as \noun{xor,
and, not}. In quantum computation, the classical bits are replaced
by qubits, which are vectors in projective Hilbert space $\mathbb{C}P^{n}$,
and the logical gates are \emph{all }the elements of the projective
unitary group $G=PU(n)$. In the real world, one must implement some
finite set of these gates, and use them to approximate the others.
Denoting by $S^{(\ell)}$ the set of $\ell$-wise products of elements
in $S\subset G$, we say that $S$ is \emph{universal }if $\left\langle S\right\rangle =\cup_{\ell\geq0}S^{(\ell)}$
is dense in $G$ (with respect to the standard bi-invariant metric
$d^{2}\left(A,B\right)=1-\frac{|\mathrm{trace}(A^{*}B)|}{2}$). This
means that any gate can be approximated with arbitrary precision as
a product of elements of $S$. The notion of Golden Gates is a much
stronger one, loosely requiring the following (see \cite{Parzanchevski2018SuperGoldenGates,Evra2018RamanujancomplexesGolden}
for precise definitions):
\begin{enumerate}
\item The covering rate of $G$ by $\left\langle S\right\rangle $ is (almost)
optimal. Namely, for every $\ell$ the set $S^{(\ell)}$ distributes
in $G$ as a perfect sphere packing (or randomly placed points) would,
up to a negligible factor.
\item Approximation: given $A\in PU(n)$ and $\varepsilon>0$, there is
an efficient algorithm to find some $A'\in B_{\varepsilon}(A)$ (the
$\varepsilon$-ball around $A$) such that $A'\in S^{(\ell)}$ with
$\ell$ (almost) minimal.
\item Compiling: given $A\in\left\langle S\right\rangle $ as a matrix,
there is an efficient algorithm to write $A$ as a word in $S$ of
the smallest possible length.
\end{enumerate}
These requirements ensure that any gate can be approximated and compiled
as an efficient circuit using the gates in $S$.

To see the connection between covering and spectral expansion, denote
by $T_{S}$ the $S$-adjacency operator on $L^{2}\left(G\right)$,
namely, $(T_{S}f)(g)=\sum_{s\in S}f(sg)$. Clearly, $T_{S}\left(\one\right)=|S|\cdot\one$,
and we denote $\lambda_{S}=\left\Vert T_{S}\big|_{\one^{\bot}}\right\Vert $,
where $\one^{\bot}=\{f\,|\,\int_{G}f\,d\mu=0\}$ and $\mu$ is the
normalized Haar measure on $G$.
\begin{thm*}[{\cite[§3]{Parzanchevski2018SuperGoldenGates}}]
Denoting by $\mu_{\varepsilon}=\mu(B_{\varepsilon}(1))$ the volume
of an $\varepsilon$-ball in $G$, the $\varepsilon$-neighborhood
of $S$ satisfies
\[
\mu\left(\bigcup\nolimits _{s\in S}B_{\varepsilon}\left(s\right)\right)\geq1-\frac{\lambda_{S}^{2}}{|S|^{2}\mu_{\varepsilon}}.
\]
\end{thm*}
Thus, as in the case of expander graphs, one aims to minimize the
nontrivial eigenvalues of an adjacency operator. It turns out that
the spectral bounds for Ramanujan graphs reappear in these settings:
\begin{thm*}[\cite{lubotzky1986hecke,lubotzky1987hecke}]
\begin{enumerate}
\item If $S\subset PU\left(2\right)$ is a symmetric set of size $k$, then
$\lambda_{S}\geq2\sqrt{k-1}$.
\item For $p\equiv1\Mod{4}$, there is an explicit symmetric set $S_{p}\subset PU\left(2\right)$
of size $k=p+1$ such that $\lambda_{S_{p}}=2\sqrt{k-1}$.
\end{enumerate}
\end{thm*}
In fact, the connection to Ramanujan graphs runs deeper than the spectral
bound. The construction of $S_{p}$, and of the $\left(p+1\right)$-regular
Ramanujan graphs in \cite{LPS88}, can be described using a single
subgroup of $PU_{2}\left(\mathbb{Q}\right)$, which acts simply-transitively
on the Bruhat-Tits tree of $PU_{2}\left(\mathbb{Q}_{p}\right)\cong PGL_{2}\left(\mathbb{Q}_{p}\right)$
(this isomorphism follows from $p\equiv1\Mod{4}$). This also solves
the compiling problem: by writing any $A\in\left\langle S_{p}\right\rangle $
in $p$-adic coordinates, one recovers its decomposition in $\left\langle S_{p}\right\rangle $
by following the (unique) path leading from $A$ to the root of the
tree (cf.\ \cite{Parzanchevski2018SuperGoldenGates}).

The proof of the spectral bound $\lambda_{S_{p}}=2\sqrt{k-1}$ uses
again the Ramanujan-Peterson conjecture (Deligne's theorem), but while
\cite{LPS88} uses the RP conjecture for automorphic representations
of weight two and arbitrary level, \cite{lubotzky1986hecke,lubotzky1987hecke}
use the conjecture for representations of level two and arbitrary
weight. To see that the gates of \cite{lubotzky1986hecke,lubotzky1987hecke}
are optimally covering (compared with random ones), one needs to bound
$\lambda_{S^{(\ell)}}$ for general $\ell$; we refer the reader to
\cite{Parzanchevski2018SuperGoldenGates} for a full account, which
addresses also the approximation problem for these gates by the Ross-Selinger
algorithm \cite{ross2015optimal}.

As Ramanujan graphs appear when studying $PU(2)$, one expects Ramanujan
complexes to appear when moving to general $PU(n)$. This is indeed
so, but the direction taken in §\ref{sec:Ramanujan-complexes}, of
replacing $\mathbb{Q}$ by $\mathbb{F}_{p}\left(t\right)$, cannot
be used anymore, since the latter does not embed in $\mathbb{R}$.
The task also becomes more complicated due to the fact that the naive
generalization of RP conjecture to $PGL_{d}$ fails, due to the appearance
of functorial lifts (cf.\ \cite{Sarnak2005NotesgeneralizedRamanujan}).
For general $n$, this is still work in progress, but for $PU(3)$
(which corresponds to quantum computation on a single \emph{qutrit}),
a complete solution exists:
\begin{thm*}[\cite{Evra2018RamanujancomplexesGolden}]
For $p\equiv1\Mod{4}$, there is an explicit Golden Gate set $S_{p}\subset PU_{3}(\mathbb{Q})$,
such that $\left\langle S_{p}\right\rangle $ acts simply transitively
on the Bruhat-Tits building of $PGL_{3}\left(\mathbb{Q}_{p}\right)$.
\end{thm*}
The compiling problem for these gates is solved by studying their
action on the two dimensional building of $PGL_{3}\left(\mathbb{Q}_{p}\right)$.
The optimal covering rate is obtained by showing that the spectral
bound $\lambda_{S_{p}}$ is the same as the maximal non-trivial adjacency
eigenvalue of a two-dimensional Ramanujan complex! Let us mention
that the proof of this bound uses Rogawski's work \cite{Rogawski1990Automorphicrepresentationsunitary},
as well as some state-of-the-art results of the Langlands program,
in particular Ngô's proof of the Fundamental Lemma, which enabled
Shin to prove the RP conjecture for cuspidal self-dual representations
of $PGL_{d}$ over CM fields \cite{Shin2011Galoisrepresentationsarising}.

\bibliographystyle{amsplain}
\bibliography{/home/parzan/Math/mybib}

\begin{bibdiv}
\begin{biblist}

\bib{ballantine2015explicit}{incollection}{
      author={Ballantine, C.},
      author={Feigon, B.},
      author={Ganapathy, R.},
      author={Kool, J.},
      author={Maurischat, K.},
      author={Wooding, A.},
       title={Explicit construction of {R}amanujan bigraphs},
        date={2015},
   booktitle={Women in numbers europe},
      series={Assoc. Women Math. Ser.},
   publisher={Springer},
       pages={1\ndash 16},
}

\bib{Chapman2019CutoffRamanujancomplexes}{article}{
      author={Chapman, Michael},
      author={Parzanchevski, Ori},
       title={Cutoff on {R}amanujan complexes and classical groups},
        date={2019},
     journal={arXiv:1901.09383},
}

\bib{Harpe1993spectrumsumgenerators}{article}{
      author={de~la Harpe, P.},
      author={Robertson, A.~G.},
      author={Valette, A.},
       title={On the spectrum of the sum of generators of a finitely generated
  group, {II}},
        date={1993},
        ISSN={0010-1354},
     journal={Colloquium Mathematicum},
      volume={65},
      number={1},
       pages={87\ndash 102},
}

\bib{Evra2018RamanujancomplexesGolden}{article}{
      author={Evra, Shai},
      author={Parzanchevski, Ori},
       title={Ramanujan complexes and {G}olden {G}ates in ${PU} (3)$},
        date={2018},
     journal={arXiv:1810.04710},
}

\bib{first2016ramanujan}{article}{
      author={First, U.A.},
       title={The {R}amanujan property for simplicial complexes},
        date={2016},
     journal={arXiv:1605.02664},
}

\bib{friedman2003relative}{article}{
      author={Friedman, J.},
       title={Relative expanders or weakly relatively {R}amanujan graphs},
        date={2003},
     journal={Duke Mathematical Journal},
      volume={118},
      number={1},
       pages={19\ndash 35},
}

\bib{golubev2017cutoff}{article}{
      author={Golubev, Konstantin},
      author={Kamber, Amitay},
       title={Cutoff on hyperbolic surfaces},
        date={2017},
     journal={arXiv:1712.10149},
}

\bib{greenberg1995spectrum}{thesis}{
      author={Greenberg, Yoseph},
       title={On the spectrum of graphs and their universal covering},
        type={Ph.D. Thesis},
        date={1995},
}

\bib{grigorchuk1999asymptotic}{article}{
      author={Grigorchuk, R.I.},
      author={{\.Z}uk, A.},
       title={On the asymptotic spectrum of random walks on infinite families
  of graphs},
        date={1999},
     journal={Random walks and discrete potential theory (Cortona, 1997),
  Sympos. Math},
      volume={39},
       pages={188\ndash 204},
}

\bib{hashimoto1989zeta}{inproceedings}{
      author={Hashimoto, K.},
       title={Zeta functions of finite graphs and representations of $p$-adic
  groups},
        date={1989},
   booktitle={Automorphic forms and geometry of arithmetic varieties},
      series={Advanced Studies in Pure Mathematics},
      volume={15},
       pages={211\ndash 280},
}

\bib{HLW06}{article}{
      author={Hoory, S.},
      author={Linial, N.},
      author={Wigderson, A.},
       title={Expander graphs and their applications},
        date={2006},
     journal={Bulletin of the American Mathematical Society},
      volume={43},
      number={4},
       pages={439\ndash 562},
}

\bib{kamber2016lp}{article}{
      author={Kamber, Amitay},
       title={${L}_p$-expander complexes},
        date={2016},
     journal={arXiv:1701.00154},
}

\bib{lafforgue2002chtoucas}{article}{
      author={Lafforgue, L.},
       title={Chtoucas de {D}rinfeld et correspondance de {L}anglands},
        date={2002},
     journal={Inventiones mathematicae},
      volume={147},
      number={1},
       pages={1\ndash 241},
}

\bib{LiRamanujanconjectureits}{article}{
      author={Li, W.C.W.},
       title={The {R}amanujan conjecture and its applications},
     journal={This volume},
}

\bib{li2004ramanujan}{article}{
      author={Li, W.C.W.},
       title={Ramanujan hypergraphs},
        date={2004},
     journal={Geometric and Functional Analysis},
      volume={14},
      number={2},
       pages={380\ndash 399},
}

\bib{Lovasz1996Randomwalksgraphs}{incollection}{
      author={Lov\'{a}sz, L.},
       title={Random walks on graphs: a survey},
        date={1996},
   booktitle={Combinatorics, {P}aul {E}rd{\"o}s is eighty, {V}ol. 2
  ({K}eszthely, 1993)},
      series={Bolyai Soc. Math. Stud.},
      volume={2},
   publisher={J\'{a}nos Bolyai Math. Soc., Budapest},
       pages={353\ndash 397},
}

\bib{Lubetzky2017RandomWalks}{article}{
      author={Lubetzky, E.},
      author={Lubotzky, A.},
      author={Parzanchevski, O.},
       title={Random walks on {R}amanujan complexes and digraphs},
        date={2019},
     journal={Journal of the European Mathematical Society, to appear},
        note={arXiv:1702.05452},
}

\bib{lubetzky2016cutoff}{article}{
      author={Lubetzky, Eyal},
      author={Peres, Yuval},
       title={Cutoff on all {R}amanujan graphs},
        date={2016},
     journal={Geometric and Functional Analysis},
      volume={26},
      number={4},
       pages={1190\ndash 1216},
}

\bib{lubetzky2011explicit}{article}{
      author={Lubetzky, Eyal},
      author={Sly, Allan},
       title={Explicit expanders with cutoff phenomena},
        date={2011},
     journal={Electronic Journal of Probability},
      volume={16},
       pages={419\ndash 435},
}

\bib{Lub10}{book}{
      author={Lubotzky, A.},
       title={Discrete groups, expanding graphs and invariant measures},
      series={Modern Birkh\"auser Classics},
   publisher={Birkh\"auser Verlag, Basel},
        date={1994},
         url={http://dx.doi.org/10.1007/978-3-0346-0332-4},
        note={With an appendix by Jonathan D. Rogawski},
}

\bib{Lub12}{article}{
      author={Lubotzky, A.},
       title={Expander graphs in pure and applied mathematics},
        date={2012},
     journal={Bull. Amer. Math. Soc},
      volume={49},
       pages={113\ndash 162},
}

\bib{Lubotzky2013}{article}{
      author={Lubotzky, A.},
       title={Ramanujan complexes and high dimensional expanders},
        date={2014},
        ISSN={0289-2316},
     journal={Japanese Journal of Mathematics},
      volume={9},
      number={2},
       pages={137\ndash 169},
         url={http://dx.doi.org/10.1007/s11537-014-1265-z},
}

\bib{lubotzky2017high}{inproceedings}{
      author={Lubotzky, A.},
       title={High dimensional expanders},
        date={2019},
   booktitle={Proc. {I}nt. {C}ong. of {M}ath. -- 2018},
      volume={1},
       pages={705\ndash 730},
        note={arXiv:1712.02526},
}

\bib{lubotzky1986hecke}{article}{
      author={Lubotzky, A.},
      author={Phillips, R.},
      author={Sarnak, P.},
       title={Hecke operators and distributing points on the sphere {I}},
        date={1986},
     journal={Comm.\ Pure.\ Appl.\ Math.},
      volume={39},
      number={1},
       pages={149\ndash 186},
}

\bib{lubotzky1987hecke}{article}{
      author={Lubotzky, A.},
      author={Phillips, R.},
      author={Sarnak, P.},
       title={Hecke operators and distributing points on {$S^2$}. {II}},
        date={1987},
     journal={Comm.\ Pure.\ Appl.\ Math.},
      volume={40},
      number={4},
       pages={401\ndash 420},
}

\bib{LPS88}{article}{
      author={Lubotzky, A.},
      author={Phillips, R.},
      author={Sarnak, P.},
       title={Ramanujan graphs},
        date={1988},
     journal={Combinatorica},
      volume={8},
      number={3},
       pages={261\ndash 277},
}

\bib{Lubotzky2005a}{article}{
      author={Lubotzky, A.},
      author={Samuels, B.},
      author={Vishne, U.},
       title={{ \!}{R}amanujan complexes of type $\tilde{A}_{d}$},
        date={2005},
     journal={Israel J. Math.},
      volume={149},
      number={1},
       pages={267\ndash 299},
}

\bib{Lubotzky2005b}{article}{
      author={Lubotzky, A.},
      author={Samuels, B.},
      author={Vishne, U.},
       title={Explicit constructions of {R}amanujan complexes of type
  $\tilde{A}_{d}$},
        date={2005},
        ISSN={0195-6698},
     journal={Eur. J. Comb.},
      volume={26},
      number={6},
       pages={965\ndash 993},
}

\bib{Lubotzky1998Noteveryuniform}{article}{
      author={Lubotzky, Alexander},
      author={Nagnibeda, Tatiana},
       title={Not every uniform tree covers {R}amanujan graphs},
        date={1998},
     journal={Journal of Combinatorial Theory, Series B},
      volume={74},
      number={2},
       pages={202\ndash 212},
}

\bib{marcus2013interlacing}{article}{
      author={Marcus, A.},
      author={Spielman, D.A.},
      author={Srivastava, N.},
       title={Interlacing families {I}: Bipartite {R}amanujan graphs of all
  degrees},
        date={2015},
     journal={Annals of Mathematics},
      volume={182},
       pages={307\ndash 325},
}

\bib{margulis1988explicit}{article}{
      author={Margulis, G.A.},
       title={Explicit group-theoretical constructions of combinatorial schemes
  and their application to the design of expanders and concentrators},
        date={1988},
     journal={Problemy Peredachi Informatsii},
      volume={24},
      number={1},
       pages={51\ndash 60},
}

\bib{morgenstern1994existence}{article}{
      author={Morgenstern, M.},
       title={Existence and explicit constructions of $q+1$ regular {R}amanujan
  graphs for every prime power $q$},
        date={1994},
     journal={Journal of Combinatorial Theory, Series B},
      volume={62},
      number={1},
       pages={44\ndash 62},
}

\bib{Parzanchevski2018SuperGoldenGates}{article}{
      author={Parzanchevski, O.},
      author={Sarnak, P.},
       title={{S}uper-{G}olden-{G}ates for {$PU(2)$}},
        date={2018},
        ISSN={0001-8708},
     journal={Advances in Mathematics},
      volume={327},
       pages={869 \ndash  901},
  url={http://www.sciencedirect.com/science/article/pii/S0001870817301640},
        note={Special volume honoring David Kazhdan},
}

\bib{Parzanchevski2018RamanujanGraphsDigraphs}{inproceedings}{
      author={Parzanchevski, Ori},
       title={Ramanujan graphs and digraphs},
        date={2019},
   booktitle={{A}nalysis and geometry on graphs and manifolds ({LMS} {L}ecture
  {N}otes), to appear},
      series={LMS Lecture Notes},
        note={arXiv:1804.08028},
}

\bib{RogawskiModularformsRamanujan}{misc}{
      author={Rogawski, J.~D.},
       title={Modular forms, the {R}amanujan conjecture and the
  {J}acquet-{L}anglands correspondence},
        note={Appendix to \cite{Lub10}},
}

\bib{Rogawski1990Automorphicrepresentationsunitary}{book}{
      author={Rogawski, J.~D.},
       title={Automorphic representations of unitary groups in three
  variables},
      series={Annals of Mathematics Studies},
   publisher={Princeton University Press, Princeton, NJ},
        date={1990},
      volume={123},
        ISBN={0-691-08586-2; 0-691-08587-0},
         url={https://doi.org/10.1515/9781400882441},
}

\bib{ross2015optimal}{article}{
      author={Ross, Neil~J},
      author={Selinger, Peter},
       title={Optimal ancilla-free {C}lifford+{V} approximation of
  z-rotations},
        date={2015},
     journal={Quantum Information \& Computation},
      volume={15},
      number={11-12},
       pages={932\ndash 950},
}

\bib{Sarnak2005NotesgeneralizedRamanujan}{incollection}{
      author={Sarnak, P.},
       title={Notes on the generalized {R}amanujan conjectures},
        date={2005},
   booktitle={Harmonic analysis, the trace formula, and {S}himura varieties},
      series={Clay Math. Proc.},
      volume={4},
       pages={659\ndash 685},
}

\bib{sarveniazi2007explicit}{article}{
      author={Sarveniazi, A.},
       title={Explicit construction of a {R}amanujan
  $\left(n_{1},n_{2},\ldots,n_{d-1}\right)$-regular hypergraph},
        date={2007},
     journal={Duke Mathematical Journal},
      volume={139},
      number={1},
       pages={141\ndash 171},
}

\bib{satake1966spherical}{inproceedings}{
      author={Satake, I.},
       title={Spherical functions and {R}amanujan conjecture},
        date={1966},
   booktitle={Proc. sympos. pure math},
      volume={9},
       pages={258\ndash 264},
}

\bib{Shin2011Galoisrepresentationsarising}{article}{
      author={Shin, Sug~Woo},
       title={Galois representations arising from some compact {S}himura
  varieties},
        date={2011},
        ISSN={0003-486X},
     journal={Annals of Mathematics},
      volume={173},
      number={3},
       pages={1645\ndash 1741},
         url={https://doi.org/10.4007/annals.2011.173.3.9},
      review={\MR{2800722}},
}

\end{biblist}
\end{bibdiv}

\noindent \begin{flushleft}
\noun{\small{}Einstein Institute of Mathematics, Hebrew University
of Jerusalem, Jerusalem 91904, Israel.}\texttt{}~\\
\texttt{\small{}alex.lubotzky@mail.huji.ac.il, parzan@math.huji.ac.il.}{\small\par}
\par\end{flushleft}
\end{document}